\documentclass[12pt]{article}
\usepackage{fancyhdr,amsmath, amssymb, psfrag, graphicx, amsfonts, layout, color,enumerate,mathrsfs}
\usepackage[numbers,sort]{natbib}

\usepackage[scaled=.95]{helvet}
\usepackage{enumitem,tikz}
\setenumerate[1]{label=\it{\textbf{\roman*)} }}

\usepackage{ifpdf}
\ifpdf 
\DeclareGraphicsExtensions{.pdf,.pdftex_t,.png,.jpg}
\else 
\DeclareGraphicsExtensions{.eps,.ps,.pst}
\fi
\usepackage[english]{babel}

\font\fontauthors=cmcsc10 scaled \magstep1
\font\footsc=cmcsc10 at 8truept
 at 8truept

\newtheorem{Th}{Theorem}

\newtheorem{Lem}[Th]{Lemma}
\newtheorem{Prop}[Th]{Proposition}
\newtheorem{Rem}[Th]{Remark}

\newcommand{\pff}{\noindent {\sc Proof.}\ }
\def\QED{\hfill$\sqcap\!\!\!\!\sqcup$}

\newtheorem{prerem}[Th]{Remark}

\newcommand{\limite}[2]{\mathop{\longrightarrow}
\limits_{\mathrm{#1}}^{\mathrm{#2}}}

\def\equivalent#1{{\mathop{\sim }\limits_{#1}}~}

\newcommand{\lpref}{\smash{\raisebox{3.5pt}{\!\!\!\begin{tabular}{c}$\hskip-4pt\scriptstyle\longleftarrow$ \\[-7pt]{\rm pref}\end{tabular}\!\!}}}
\newcommand{\petitlpref}{\smash{\raisebox{2.5pt}{\!\!\!\begin{tabular}{c}$\hskip-2pt\scriptscriptstyle\longleftarrow$ \\[-9pt]{$\scriptstyle\hskip 1pt\rm pref$}\end{tabular}\!\!}}}

\def\Var{\mathop{\rm Var}\nolimits}
\def\stirling2 #1#2{\left\{\begin{matrix} #1\\#2\end{matrix}\right\}}

\newcommand\1{\leavevmode\hbox{\rm \small1\kern-0.35em\normalsize1}}
\newcommand\ind[1]{\1_{\{#1\}}}
\newcommand{\petito}[1]{o\mathopen{}\left(#1\right)}
\newcommand{\grandO}[1]{\mathcal{O}\mathopen{}\left(#1\right)}

\newcommand{\N}{{{\rm I}\!{\rm N}}}
\newcommand{\Q}{{\rm Q}\kern-.65em{}^{{}_/}}
\newcommand{\PP}{\mathbb{P}}
\def\1{{\bf 1}}

\def\g#1{\mathbb #1}
\def\rond#1{\mathcal #1}

\def\N{\g N}

\usepackage[draft]{fixme}
\fxusetheme{color}
\FXRegisterAuthor{p}{apa}{P}
\FXRegisterAuthor{b}{ape}{B}
\FXRegisterAuthor{f}{af}{F}
\FXRegisterAuthor{n}{an}{N}
\begin{document}

\begin{center}
\LARGE{\bf {Uncommon Suffix Tries}}
\end{center}

\begin{center}
{\fontauthors
Peggy C\'enac}\footnote{
Universit\'e de Bourgogne,
Institut de Math\'ematiques de Bourgogne,
IMB UMR 5584 CNRS,
9 rue Alain Savary - BP 47870, 21078 DIJON CEDEX, France.
}

\medskip
{\fontauthors
Brigitte Chauvin}\footnote{
Universit\'e de Versailles-St-Quentin,
Laboratoire de Math\'ematiques de Versailles,
CNRS, UMR 8100,
45, avenue des Etats-Unis, 78035 Versailles CEDEX, France.
}

\medskip
{\fontauthors
Fr\'ed\'eric Paccaut}\footnote{
LAMFA,
CNRS, UMR 6140,
Universit\'e  de Picardie Jules Verne,
33, rue Saint-Leu, 80039 Amiens, France.
}

\medskip
{\fontauthors
Nicolas Pouyanne}\footnote{
Universit\'e de Versailles-St-Quentin,
Laboratoire de Math\'ematiques de Versailles,
CNRS, UMR 8100,
45, avenue des Etats-Unis, 78035 Versailles CEDEX, France.
}

\vskip 20pt
{\it December 14th 2011}
\end{center}

\vskip 20pt
{\small
{\bf Abstract}

\medskip

Common assumptions on the source producing the words inserted in a suffix trie with $n$ leaves lead to a $\log n$ height and saturation level. We provide an example of a suffix trie whose height increases faster than a power of $n$ and another one whose saturation level is negligible with respect to $\log n$. Both are built from VLMC (Variable Length Markov Chain) probabilistic sources and are easily extended to families of tries having the same properties. The first example corresponds to a ``logarithmic infinite comb'' and enjoys a non uniform polynomial mixing. The second one corresponds to a ``factorial infinite comb'' for which mixing is uniform and exponential.
\vskip 10pt
{\it MSC 2010}: 60J05, 37E05.

\vskip 5pt
{\it Keywords}: variable length Markov chain, probabilistic source, mixing properties, suffix trie 
}


\section{Introduction}
\label{sec:Intro}

\emph{Trie} (abbreviation of re{\it trie}val) is a natural data structure, efficient for searching words in a given set and used in many algorithms as data compression, spell checking or IP addresses lookup. A \emph{trie} is a digital tree in which
words are inserted in external nodes. The trie process grows up by successively inserting words according to their prefixes. A precise definition will be given in Section~\ref{ssec:height}.

As soon as a set of
words is given, the way they are inserted in the trie is deterministic. Nevertheless, a trie becomes random when the words are randomly drawn: each word is produced by a probabilistic source and $n$ words are chosen (usually independently) to be inserted in a trie. A \emph{suffix trie} is a trie built on the suffixes of \emph{one} infinite word. The randomness then comes from the source producing such an infinite word and the successive words inserted in the tree are far from being independent, they are strongly correlated. 

As a principal application of suffix tries one can cite the lossless compression algorithm Lempel-Ziv 77 (LZ77). The first results on the average size of suffix tries when the infinite word is given by a symmetrical memoryless source are due to Blumer et al.\@ \cite{Bl89} and those on the height of the tree to Devroye \cite{Devroye2}. Using analytic combinatorics, Fayolle \cite{Fayolle} has obtained the average size and total path length of the tree for a binary word issued from a memoryless source (with some restriction on the probability of each letter). 

Here we are interested in the height $H_n$ and the saturation level $\ell_n$ of a suffix trie $\rond T_n$ containing the first $n$ suffixes of an infinite word produced by a source associated with a so-called Variable Length Markov Chain (VLMC) (see Rissanen \cite{Ris} for the seminal work, Galves-L\"ocherbach \cite{GalLoc} for an overview,  and  \cite{ccpp} for a probabilistic frame).  One deals with a particular VLMC source associated with an infinite comb, described hereafter. This particular model has the double advantage to go beyond the cases of memoryless or Markov sources and to provide concrete computable properties. The analysis of the height and the saturation level is usually motivated by optimization of the memory cost. Height is clearly relevant to this point; saturation level is algorithmically relevant as well because internal nodes below the saturation level are often replaced by a less expansive table.

All the tries or suffix tries considered so far in the literature have a height and a saturation level both growing logarithmically with the number of words inserted, to the best of our knowledge.
For plain tries, when the inserted words are independent, the results due to Pittel \cite{Pittel} rely on  two assumptions on the source producing the words: first, the source is uniformly mixing, second, the probability of any word decays exponentially with its length.
Let us also mention the general analysis of tries by Cl\'ement-Flajolet-Vall\'ee \cite{ClementFlajoletVallee} for dynamical sources.
For suffix tries, Szpankowski \cite{Szpan} obtains the same result, with a weaker mixing assumption (still uniform though) and the same hypothesis on the measure of the words.

Our aim is to exhibit two cases when these behaviours are no longer the same. The first example is the ``\emph{logarithmic comb}'', for which we show that the mixing is slow in some sense, namely non uniformly polynomial (see Section \ref{ssec:mixinglogarithmiccomb} for a precise statement) and the measure of some increasing sequence of words decays polynomially. We prove in Theorem~\ref{hauteur} that the height of this trie is larger than a power of $n$ (when $n$ is the number of inserted suffixes in the tree). The second example is the ``\emph{factorial comb}'', which has a uniformly exponential mixing, thus fulfilling the mixing hypothesis of Szpankowski  \cite{Szpan}, but the measure of some increasing sequence of words decays faster than any exponential. In this case we prove in Theorem~\ref{saturation} that the saturation level is negligible with respect to $\log n$. We prove more precisely that, almost surely, 
$\ell _n\in o\left(\frac{\log n}{(\log\log n)^\delta}\right)$, for any $\delta >1$.

\medskip

The paper is organised as follows. In Section \ref{sec:source}, we define a VLMC source associated with an infinite comb. In Section, \ref{sec:mixing} we give results on the mixing properties of these sources by explicitely computing the suitable generating functions in terms of the source data. In Section \ref{sec:suffixtrie}, the associated suffix tries are built, and the two uncommon behaviours are stated and shown. The methods are based on two key tools concerning pattern return time: a duality property and the computation of generating functions. The relation between the mixing of the source and the asymptotic behaviour of the trie is highlighted by the proof of Proposition~\ref{Phi2}.

\section{Infinite combs as sources}
\label{sec:source}

In this section, a VLMC probabilistic source associated with an infinite comb is defined. Moreover, we introduce the two examples given in introduction: the logarithmic and the factorial combs. We begin with the definition of a general variable length Markov Chain associated with a probabilized infinite comb. 

\begin{figure}[h]
\begin{center}
\includegraphics[height=70 truemm]{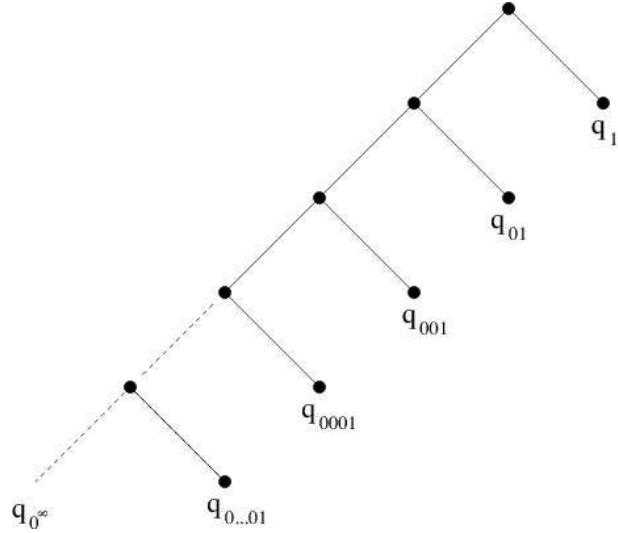}
\end{center}
\caption{\label{peigne}An infinite comb}
\end{figure}

The following presentation comes from \cite{ccpp}. Let $\rond A$ be the alphabet $\{0,1\}$ and $\rond L = \rond A^{-\N}$ be the set of left-infinite words. Consider the binary tree (represented in Figure \ref{peigne}) whose finite leaves are the words $1, 01, \dots, 0^k1, \dots $ and with an infinite leaf $0^{\infty}$ as well. Each leaf is labelled with a Bernoulli distribution, respectively  denoted by $q_{0^k1}, k\geqslant 0$ and $q_{0^\infty}$. This probabilized tree is called \emph{the infinite comb}.

The VLMC (Variable Length Markov Chain) associated with an infinite comb is the $\rond L$-valued Markov chain $(V_n)_{n\geqslant 0}$ defined by the transitions
$$\g P(V_{n+1} = V_n\alpha | V_n)= q_{\petitlpref (V_n)}(\alpha )
$$
where $\alpha\in\rond A$ is any letter and $\lpref(V_n)$ denotes the first suffix of $V_n$ (reading from right to left) appearing as a leaf of the infinite comb. For instance, if $V_n = \dots 1000$, then $\lpref(V_n) = 0001$. 
Notice that the VLMC is entirely determined by the data $q_{0^{\infty}}, q_{0^k1}, k\geqslant 0$.
From now on, denote $c_0 = 1$ and for $n\geqslant 1$,
$$
c_n:= \prod_{k=0}^{n-1} q_{0^{k}1}(0).
$$

\medskip

It is proved in \cite{ccpp} that in the irreducible case \emph{i.e.}\@ when $q_{0^{\infty}}(0)\not= 1$, there exists a unique stationary probability measure $\pi$ on $\rond L$ for $(V_n)_n$ if and only if the series $\sum c_n$ converges. 
From now on, we assume that this condition is fulfilled and we call 
\begin{equation}
\label{def-S}
S(x):=\sum_{n\geqslant 0} c_n x^n
\end{equation}
its generating function so that $S(1) = \sum_{n\geqslant 0} c_n$. For any finite word $w$, we denote $\pi(w):=\pi(\rond Lw)$. Computations performed in \cite{ccpp} show that for any $n\geqslant 0$,
  \begin{equation}
 \label{probas}
 \pi(10^n)=\frac{c_n}{S(1)} \quad \mbox{and}\quad \pi(0^n)=\frac{\sum_{k \geqslant n}c_k}{S(1)}. \end{equation}
Notice that, by stationarity $\pi(0^n)=\pi(0^{n+1})+\pi(10^n)$ and by disjointness of events, $\pi(0^n)=\pi(0^{n+1})+ \pi(0^n1)$ for all $n\geqslant 1$ so that
\begin{equation}
\label{sym}
\pi(10^n) =  \pi(0^n1) .
\end{equation}
If $U_n$ denotes the final letter of $V_n$, the random sequence $W=U_0U_1U_2\ldots$ is a right-infinite random word. We define in this way a probabilistic source in the sense of information theory \emph{i.e.}\@ a mechanism that produces random words. This VLMC probabilistic source is characterized by: 
$$
p_w:=\g P(W \hbox{ has } w \hbox{ as a prefix }) = \pi(w),
$$
for every finite word $w$. Both particular suffix tries the article deals with are built from such sources, defined by the following data.

\medskip

\subsubsection*{Example 1: the logarithmic comb}

The logarithmic comb is defined by $c_0=1$ and for $n\geqslant 1$, 
\[c_n = \frac{1}{n(n+1)(n+2)(n+3)}. \]
The corresponding conditional probabilities on the leaves of the tree are
  \[q_1(0)=\frac{1}{24}\quad \mbox{and for }n\geqslant 1, \quad q_{0^n1}(0)=1-\frac{4}{n+4}.\]
The expression of $c_n$ was chosen to make the computations as simple as possible and also because the square-integrability of the waiting time of some pattern will be needed (see end of Section~\ref{ssec:generating}), guaranteed by
\[\sum_{n\geqslant 0} n^2c_n<+\infty.\]

\subsubsection*{Example 2: the factorial comb}

\medskip

The conditional probabilities on the leaves are defined by
$$
q_{0^n1}(0)= \frac 1{n+2} \quad \mbox{ for }n\geqslant 0,
$$
so that
$$
c_n = \frac 1{(n+1)!}.
$$

\section{Mixing properties of infinite combs}
\label{sec:mixing}

In this section, we first precise what we mean by mixing properties of a random sequence. We refer to Doukhan  \cite{Doukhan}, especially for the notion of $\psi$-mixing defined in that book. We state in Proposition \ref{melange} a general result that provides the mixing coefficient for an infinite comb defined by $(c_n)_{n\geqslant 0}$ or equivalently by its generating function~$S$. This result is then applied to our two examples. The mixing of the logarithmic comb is polynomial but not uniform, it is a very weak mixing; the mixing of the factorial comb is uniform and exponential, it is a very strong mixing. Notice that mixing properties of some infinite combs have already been investigated by Isola~\cite{Isola}, although with a slight different language.

\subsection{Mixing properties of general infinite combs}
\label{ssec:mixingcomb}

For a stationary sequence $(U_n)_{n\geqslant 0}$ with stationary measure $\pi$, we want to measure by means of a suitable coefficient the independence between two words $A$ and $B$ separated by $n$ letters. The sequence is said to be ``mixing'' when this coefficient vanishes when $n$ goes to $+\infty$. Among all types of mixing, we focus on one of the strongest type: \emph{$\psi$-mixing}. More precisely, for $0\leqslant m\leqslant +\infty$, denote by ${\cal F}_0^m$ the $\sigma$-algebra generated by $\{U_k , 0\leqslant k \leqslant m\}$ and introduce for $A\in{\cal F}_0^m$ and $B\in{\cal F}_0^{\infty}$ the mixing coefficient
\begin{eqnarray}
\psi(n,A,B)&:=&\frac{\pi(A\cap T^{-(m+1)-n}B)-\pi(A)\pi(B)}{\pi(A)\pi(B)}\nonumber \\
&& \nonumber \\
&=&\frac{\sum_{|w|=n}\pi(AwB)-\pi(A)\pi(B)}{\pi(A)\pi(B)},\label{mixingCoeff}
\end{eqnarray}
where $T$ is the shift map and where the sum runs over the finite words $w$ with length $|w|=n$.

A sequence $(U_n)_{n\geqslant 0}$ is called \emph{$\psi$-mixing} whenever
\[\lim_{n\to\infty}Ê\ \sup_{m\geqslant 0,A\in{\cal F}_0^m, B\in{\cal F}_0^{\infty}}|\psi(n,A,B)|=0.\]
In this definition, the convergence to zero is uniform over all words $A$ and~$B$. This is not going to be the case in our first example.
As in Isola \cite{Isola}, we widely use the renewal properties of infinite combs (see Lemma~\ref{lemme_un}) but more detailed results are needed, in particular we investigate the lack of uniformity for the logarithmic comb. 

\subsubsection*{Notations and Generating functions}

\medskip

$\bullet$ For a comb, recall that $S$ is the generating function of the nonincreasing sequence $(c_n)_{n\geqslant 0}$ defined by (\ref{def-S}).
 
$\bullet$ Set $\rho_0=0$ and for $n\geqslant 1$, 
\[\rho_n:=c_{n-1}-c_{n},\]
with generating function
$$
P(x) := \sum_{n\geqslant 1} \rho_n x^n.
$$
$\bullet$ Define the sequence $(u_n)_{n\geqslant 0}$ by $u_0=1$ and for $n\geqslant 1$,
\begin{equation}\label{un}
u_n:=\frac{\pi(U_0=1,U_{n}=1)}{\pi(1)}=\frac{1}{\pi(1)}\sum_{|w|=n-1}\pi(1w1),
\end{equation}
and let
\[U(x):=\sum_{n\geqslant 0}u_nx^n\]
denote its generating function. Hereunder is stated a key lemma that will be widely used in Proposition~\ref{melange}. In some sense, this kind of relation (sometimes called Renewal Equation) reflects the renewal properties of the infinite comb. 
\begin{Lem}
\label{lemme_un}
The sequences $(u_n)_{n\geqslant 0}$ and $(\rho_n)_{n\geqslant 0}$ are connected by the relations:
\[\forall n\geqslant 1,~~u_n=\rho_n+u_1\rho_{n-1}+\ldots+u_{n-1}\rho_1\]
and (consequently)
\[U(x) = \sum_{n=0}^{\infty}u_nx^n=\frac{1}{1-P(x)}=\frac{1}{(1-x)S(x)}.\]
\end{Lem}
\pff
For a finite word $w=\alpha_1\ldots\alpha_m$ such that $w \neq 0^m$, let $l(w)$ denote the position of the last $1$ in $w$, that is $l(w):=\max\{1\le i\le m,~\alpha_i=1\}$. Then, the sum in the expression (\ref{un}) of $u_n$ can be decomposed as follows:
\[\sum_{|w|=n-1}\pi(1w1)=\pi(10^{n-1}1)+\sum_{i=1}^{n-1}\sum_{{|w|=n-1}\atop{l(w)=i}}\pi(1w1).\]
Now, by disjoint union $\pi(10^{n-1})=\pi(10^{n-1}1)+\pi(10^{n})$, so that
\[\pi(10^{n-1}1)=\pi(1)(c_{n-1}-c_{n})=\pi(1)\rho_n.\] 
In the same way, for $w=\alpha_1\ldots\alpha_{n-1}$, if $l(w)=i$ then $\pi(1w1)=\pi(1\alpha_1\ldots\alpha_{i-1}1)\rho_{n-i}$, so that
\begin{eqnarray*}
u_n & = & \rho_n+\sum_{i=1}^{n-1}\rho_{n-i}\frac{1}{\pi(1)}\sum_{|w|=i-1}\pi(1w1) \\
 & = & \rho_n+\sum_{i=1}^{n-1}\rho_{n-i}u_{i},
\end{eqnarray*}
which leads to $U(x)=(1-P(x))^{-1}$ by summation.
\QED

\subsubsection*{Mixing coefficients}
The mixing coefficients $\psi(n,A,B)$ are expressed as the $n$-th coefficient in the series expansion of an analytic function $M^{A,B}$ which is given in terms of $S$ and~$U$. The notation $[x^n] A(x)$ means the  coefficient of $x^n$ in the power expansion of $A(x)$ at the origin.
Denote the remainders associated with the series $S(x)$ by
$$
r_n : = \sum_{k\geqslant n} c_k, \quad R_n(x):=\sum_{k\geqslant n}c_k x^k
$$
and for $a\geqslant 0$, define the ``shifted'' generating function
\begin{equation}
\label{defPa}
P_a(x):=\frac 1{c_a}\sum _{n\geqslant 1}\rho _{a+n}x^n=x+\frac{x-1}{c_ax^a}R_{a+1}(x).
\end{equation}
\begin{Prop}\label{melange}
For any finite word $A$ and any word $B$, the identity
\[\psi(n,A,B)=[x^{n+1}] M^{A,B}(x)\]
holds for the generating functions $M^{A,B}$ respectively defined by:
\begin{enumerate}
\item if $A=A'1$ and $B=1B'$ where $A'$ and $B'$ are any finite words, then $$M^{A,B}(x) = M(x) := \displaystyle\frac{S(x)-S(1)}{(x-1)S(x)};$$
\item if $A=A'10^a$ and $B=0^b1B'$ where $A'$ and $B'$ are any finite words and $a+b\geqslant 1$, then 
\[M^{A,B}(x)
:= S(1)\displaystyle\frac{c_{a+b}}{c_ac_b}P_{a+b}(x) + U(x)\left[ S(1)P_a(x)P_b(x)-S(x)\right];\]
\item if $A=0^a$ and $B=0^b$ with $a,b\geqslant 1$, then
\[M^{A,B}(x)
:= S(1)\displaystyle\frac{1}{r_ar_b}\sum _{n\geqslant 1}r_{a+b+n-1}x^n + \displaystyle U(x)
\left[ \frac{S(1)R_a(x)R_b(x)}{r_ar_bx^{a+b-2}}-S(x)\right];
\]
\item if $A=A'10^a$ and $B=0^b$ where $A'$ is any finite words and $a,b\geqslant 0$, then
$$
M^{A,B}(x)
:= S(1)\displaystyle\frac{1}{c_ar_bx^{a+b-1}}R_{a+b}(x)+  \displaystyle U(x)
\left[ \frac{S(1)P_a(x)R_b(x)}{c_ar_bx^{b-1}}-S(x)\right];
$$
\item if $A=0^a$ and $B=0^b1B'$ where $B'$ is any finite words and $a,b\geqslant 0$, then 
  \[
  M^{A,B}(x)
  := S(1)\displaystyle\frac{1}{r_ac_bx^{a+b-1}}R_{a+b}(x)+  \displaystyle U(x)
\left[ \frac{S(1)R_a(x)P_b(x)}{r_ac_bx^{a-1}}-S(x)\right].
  \]
\end{enumerate}
\end{Prop}
\begin{Rem}
It is worth noticing that the asymptotics of $\psi(n,A,B)$ may not be uniform in all words $A$ and $B$. We  call this kind of system non-uniformly $\psi$-mixing. It may happen that $\psi(n,A,B)$ goes to zero for any fixed $A$ and $B$, but (for example, in case \textbf{\itshape iii)}) the larger $a$ or $b$, the slower the convergence, preventing it from being uniform.
\end{Rem}
\pff
The following identity has been established in \cite{ccpp} (see formula (17) in that paper) and will be used many times in the sequel. For any two finite words $w$ and $w'$,
\begin{equation}
\label{renouvellement}
\pi(w1w')\pi(1)=\pi(w1)\pi(1w').
\end{equation}
\begin{enumerate}
\item  
 If $A=A'1$ and $B=1B'$, then (\ref{renouvellement}) yields 
\[\pi(AwB)=\pi(A'1w1B')=\frac{\pi(A'1)}{\pi(1)}\pi(1w1B')=S(1)\pi(A)\pi(B)\frac{\pi(1w1)}{\pi(1)}.\]
So  
$$
\psi(n,A,B) = S(1) u_{n+1} - 1
$$
and by Lemma \ref{lemme_un}, the result follows.

\item
Let $A=A'10^a$ and $B=0^b1B'$ with $a,b\geqslant 0$ and $a+b\neq 0$.
To begin with,
\[\pi(AwB)=\frac{1}{\pi(1)}\pi(A'1)\pi(10^aw0^b1B')=\frac{1}{\pi(1)^2}\pi(A'1)\pi(10^aw0^b1)\pi(1B').\]
Furthermore, $\pi(A)=c_a\pi(A'1)$ and  by (\ref{sym}), $\pi(0^b1)=\pi(10^b)$, so it comes 
 \[\pi(B)=\frac{1}{\pi(1)}\pi(0^b1)\pi(1B')=\frac{\pi(10^b)}{\pi(1)}\pi(1B')=c_b\pi(1B').\]
 Therefore,
 \[\pi(AwB)=\frac{\pi(A)\pi(B)}{c_ac_b\pi(1)^2}\pi(10^aw0^b1).\]
Using $\pi(1)S(1)=1$, this proves
\[\psi(n,A,B)=S(1)\frac{v_{n}^{a,b}}{c_ac_b}-1\]
where
  \[v_n^{a,b}:=\frac{1}{\pi(1)}\sum_{|w|=n-1}\pi(10^aw0^b1).\] 
As in the proof of the previous lemma, if $w=\alpha_1\ldots\alpha_{m}$ is any finite word different from $0^m$, we call $f(w):=\min\{1\le i \le m, \alpha_{i}=1\}$ the first place where  $1$ can be seen in $w$ and recall that $l(w)$ denotes the last place where  $1$ can be seen in $w$. One has
\[\sum_{|w|=n-1}\pi(10^aw0^b1) = \pi(10^{a+n-1+b}1)+\sum_{1\le i\le j\le n-1}\sum_{{{|w|=n-1}\atop{f(w)=i,l(w)=j}}}\pi(10^aw0^b1).\]
If $i=j$ then $w$ is the word $0^{i-1}10^{n-i-1}$, else $w$ is of the form $0^{i-1}1w'10^{n-1-j}$, with $|w'|=j-i-1$. Hence, the previous sum can be rewritten as
\begin{eqnarray*}
  \sum_{|w|=n-1}\pi(10^aw0^b1) & = & \pi(1)\rho_{a+b+n}+\pi(1)\sum_{i=1}^{n-1}\rho_{a+i}\rho_{n-i+b}\\
 & + &\sum_{1\leqslant i<j\leqslant n-1}\sum_{{{w}\atop{|w|=j-i-1}}}\pi(10^{a+i-1}1w10^{n-1-j+b}1).
\end{eqnarray*}
Equation (\ref{renouvellement}) shows
\begin{eqnarray*}
\pi(10^{a+i-1}1w10^{n-1-j+b}1) & = & \frac{\pi(10^{a+i-1}1)}{\pi(1)}\frac{\pi(1w1)}{\pi(1)}\pi(10^{n-1-j+b}1) \\
 & = & \rho_{a+i}\rho_{n-j+b}\pi(1w1).
\end{eqnarray*}
This implies:
\[v_n^{a,b}=\rho_{a+b+n}+\sum_{i=1}^{n-1}\rho_{a+i}\rho_{n-i+b}
   +  \sum_{1\le i<j\le n-1}\rho_{a+i}\rho_{n-j+b}\sum_{w, |w|=j-i+1}\frac{\pi(1w1)}{\pi(1)}.\]
Recalling that $u_0=1$, one gets
 \[v_n^{a,b}=\rho_{a+b+n}+\sum_{1\le i\le j\le n-1}\rho_{a+i}\rho_{n-j+b}u_{j-i}\]
which gives the result  \textbf{\itshape ii)} with Lemma \ref{lemme_un}.

\item
Let $A=0^a$ and $B=0^b$ with $a,b\geqslant 1$.
Set 
 \[v_n^{a,b}:=\frac{1}{\pi(1)}\sum_{|w|=n-1}\pi(0^aw0^b).\]
First, recall that, due to (\ref{probas}), $\pi(A)=\pi(1)r_a$ and $\pi(B)=\pi(1)r_b$. Consequently,
 \[\psi(n,A,B)=\frac{\pi(1)v_{n+1}^{a,b}-\pi(A)\pi(B)}{\pi(A)\pi(B)}=S(1)\frac{v_{n+1}^{a,b}}{r_ar_b}-1.\]
Let $w$ be a finite word with $|w|=n-1$. If $w=0^{n-1}$, then 
\[\pi(AwB)=\pi(0^{a+n-1+b})=\pi(1)r_{a+b+n-1}.\] 
If not, let $f(w)$ denote as before the first position of $1$ in $w$ and $l(w)$ the last one in $w$. If $f(w)=l(w)$, then
 \[\begin{array}{rcl}
 \pi(AwB) & = & \pi(0^{a+f(w)-1}10^{n-1-f(w)+b}) \\
  & = & \frac{1}{\pi(1)}\pi(0^{a+f(w)-1}1)\pi(10^{n-1-f(w)+b})=\pi(1)c_{a+f(w)-1}c_{n-1-f(w)+b}.
  \end{array}\]
 If $f(w)<l(w)$, then writing $w=w_1 \ldots w_{n-1}$, 
 \[\begin{array}{rcl}
 \pi(AwB) & = & \pi(0^{a+f(w)-1}1w_{f(w)+1}\ldots w_{l(w)-1}10^{n-1-l(w)}) \\
   & = & \frac{1}{\pi(1)^2}\pi(0^{a+f(w)-1}1)\pi(1w_{f(w)+1}\ldots w_{l(w)-1}1)\pi(10^{n-1-l(w)+b}).
   \end{array}\]
 Summing yields
 \[\begin{array}{rcl}
 v_n^{a,b} & = & r_{a+b+n-1}+\displaystyle\sum_{i=1}^{n-1}c_{a+i-1}c_{n-1+b-i}+\sum_{{{i,j=1}\atop{i<j}}}^{n-1}\sum_{{{w,}\atop{|w|=j-i-1}}}c_{a+i-1}\frac{\pi(1w1)}{\pi(1)}c_{n-1+b-j} \\
 & = & r_{a+b+n-1}+\displaystyle\sum_{1\le i\le j\le n-1}c_{a+i-1}c_{n-1+b-j}u_{j-i},
 \end{array}\]
which gives the desired result. The last two items, left to the reader, follow the same guidelines.
\QED
\end{enumerate}

\subsection{Mixing of the logarithmic infinite comb}
\label{ssec:mixinglogarithmiccomb}
Consider the first example in Section~\ref{sec:source}, that is the probabilized infinite comb defined by $c_0 = 1$ and  for any $n\geqslant 1$ by
$$
c_n= \frac{1}{n(n+1)(n+2)(n+3)}.
$$
When $|x|<1$, the series $S(x)$ writes as follows
\begin{equation}
\label{defSpeigneLog}
S(x)=\frac{47}{36}-\frac{5}{12x}+\frac{1}{6x^2}+\frac{(1-x)^3\log(1-x)}{6x^3}
\end{equation}
and
 \[S(1)=\frac{19}{18}.\]
With Proposition~\ref{melange}, the asymptotics of the mixing coefficient comes from singularity analysis of the generating functions $M^{A,B}$. 

\begin{Prop}
\label{prop-3-4}
The VLMC defined by the logarithmic infinite comb has a non-uniform polynomial mixing of the following form:
 for any finite words $A$ and $B$, there exists a positive constant $C_{A,B}$ such that for any $n\geqslant 1$,
$$
|\psi (n,A,B)|\leqslant\frac {C_{A,B}}{n^3}.
$$
\end{Prop}
\begin{Rem}
The $C_{A,B}$ cannot be bounded above by some constant that does not depend on $A$ and $B$, as can be seen hereunder in the proof.
Indeed, we show that if $a$ and $b$ are positive integers,
$$
\psi(n,0^a,0^b)\sim\frac 13\left( \frac{S(1)}{r_ar_b}-\frac 1{r_a}-\frac 1{r_b}+\frac 1{S(1)}\right)\frac 1{n^3}
$$
as $n$ goes to infinity. In particular, $\psi(n,0,0^n)$ tends to the positive constant $\frac{13}{6}$.
\end{Rem}

\noindent {\sc Proof of Proposition~\ref{prop-3-4}.}\

For any finite words $A$ and $B$ in case \textbf{\itshape i)} of Proposition~\ref{melange}, one deals with $U(x)=\left((1-x)S(x)\right)^{-1}$ which has $1$ as a unique dominant singularity. Indeed, $1$ is the unique dominant singularity of $S$, so that the dominant singularities of $U$ are $1$ or zeroes of $S$ contained in the closed unit disc. But $S$ does not vanish on the closed unit disc, because for any $z$ such that $ |z| \leqslant 1$, 
$$|S(z)| \geqslant 1- \sum_{n\geqslant 1} \frac{1}{n(n+1)(n+2)(n+3)} = 1-(S(1)-1)=\frac{17}{18}.$$

\medskip

Since 
\[M(x) =\frac{S(x)-S(1)}{(x-1)S(x)}= S(1)U(x) - \frac 1{1-x},\]  
the unique dominant singularity of $M$ is $1$, and when $x$ tends to $1$ in the unit disc, (\ref{defSpeigneLog}) leads to
\[
M(x) = A(x)-\frac{1}{6S(1)}(1-x)^2\log(1-x)+\grandO{(1-x)^3\log(1-x)}
\]
where $A(x)$ is a polynomial of degree $2$. Using the classical transfer theorem (see Flajolet and Sedgewick \cite[section VI]{FS}) based on the analysis of the singularities of $M$, we get
 \[\psi(n-1,w1,1w')=  [x^n] M(x)=\frac{1}{3S(1)}\frac{1}{n^3}+\petito{\frac{1}{n^3}}.\]
 
The cases  \textbf{\itshape{ii)}}, \textbf{\itshape{iii)}}, \textbf{\itshape{iv)}} and \textbf{\itshape{v)}} of Proposition~\ref{melange} are of the same kind, and we completely deal with case  \textbf{\itshape{iii)}}.
 
\vskip 5pt
Case  \textbf{\itshape{iii)}}: words of the form $A=0^a$ and $B=0^b$, $a,b\geqslant 1$.
As shown in Proposition~\ref{melange}, one has to compute the asymptotics of the $n$-th coefficient of the Taylor
series of the function
\begin{equation}
\label{melangeLog3}
M^{a,b}(x):=S(1)\displaystyle\frac{1}{r_ar_b}\sum _{n\geqslant 1}r_{a+b+n-1}x^n + \displaystyle U(x)\left[ \frac{S(1)R_a(x)R_b(x)}{r_ar_bx^{a+b-2}}-S(x)\right].
\end{equation}
The contribution of the left-hand term of this sum is directly given by the asymptotics of the remainder
$$
r_n=\sum _{k\geqslant n}c_k=\frac{1}{3n(n+1)(n+2)}=\frac{1}{3n^3}+\mathcal{O}\left( \frac {1}{n^4}\right) .
$$
By means of singularity analysis, we deal with the right-hand term
$$
N^{a,b}(x):=\displaystyle U(x)\left[ \frac{S(1)R_a(x)R_b(x)}{r_ar_bx^{a+b-2}}-S(x)\right].
$$
Since $1$ is the only dominant singularity of $S$ and $U$ and consequently of any $R_a$, it suffices to compute
an expansion of $N^{a,b}(x)$ at $x=1$.
It follows from~(\ref{defSpeigneLog}) that
$U$, $S$ and $R_a$ admit expansions near $1$ of the forms
\begin{align*}
U(x)
&=\frac{1}{S(1)(1-x)}+\mbox{polynomial}+\frac{1}{6S(1)^2}(1-x)^2\log (1-x)+\mathcal{O}(1-x)^2, \\
S(x)&=\mbox{polynomial}+\frac{1}{6}(1-x)^3\log (1-x)+\mathcal{O}(1-x)^3,
\end{align*}
and
\begin{align*}
R_a(x)
&=\mbox{polynomial}+\frac{1}{6}(1-x)^3\log (1-x)+\mathcal{O}(1-x)^3. \quad
\end{align*}
Consequently,
$$
N^{a,b}(x)=\frac 16\left( \frac 1{r_a}+\frac 1{r_b}-\frac 1{S(1)}\right)(1-x)^2\log (1-x)+\mathcal{O}(1-x)^2
$$
in a neighbourhood of $1$ in the unit disc so that, by singularity analysis,
$$
[x^n]N^{a,b}(x)=-\frac 13\left( \frac 1{r_a}+\frac 1{r_b}-\frac 1{S(1)}\right)\frac 1{n^3}+o\left(\frac 1{n^3}\right).
$$
Consequently (\ref{melangeLog3}) leads to
\[
\psi(n-1,0^a,0^b)=  [x^n] M^{a,b}(x)
\sim\frac 13\left( \frac{S(1)}{r_ar_b}-\frac 1{r_a}-\frac 1{r_b}+\frac 1{S(1)}\right)\frac 1{n^3}
\]
as $n$ tends to infinity, showing the mixing inequality and the non uniformity.

The remaining cases \textbf{\itshape{ii)}}, \textbf{\itshape{iv)}} and \textbf{\itshape{v)}} are of the same flavour.
\QED
\subsection{Mixing of the factorial infinite comb}
\label{ssec:mixingfactccomb}

Consider now the second Example in Section~\ref{sec:source}, that is the probabilized infinite comb defined by
$$
\forall n\in\g N,~c_n=\frac 1{(n+1)!}.
$$
With previous notations, one gets
$$
S(x)=\frac {e^x-1}x
{\rm ~~and~~}
U(x)=\frac {x}{(1-x)(e^x-1)}.
$$

\begin{Prop}
The VLMC defined by the factorial infinite comb has a uniform exponential mixing of the following form:
there exists a positive constant $C$ such that for any $n\geqslant 1$ and for any finite words $A$ and $B$, 
$$
|\psi (n,A,B)|\leqslant\frac C{(2\pi )^n}.
$$
\end{Prop}
\pff
\begin{enumerate}
\item First case of mixing in Proposition~\ref{melange}: $A=A'1$ and $B=1B'$.

Because of Proposition~\ref{melange}, the proof consists in computing the asymptotics of $[x^n] M(x)$. We make use of singularity analysis. The dominant singularities of
$$
M(x)=\frac{S(x)-S(1)}{(x-1)S(x)}
$$
are readily seen to be $2i\pi$ and $-2i\pi$, and
$$
M(x)\ \equivalent{2i\pi}\frac{1-e}{1-2i\pi}\cdot\frac{1}{1-\frac{z}{2i\pi}}.
$$
The behaviour of $M$ in a neighbourhood of $-2i\pi$ is obtained by complex conjugacy.
Singularity analysis via  transfer theorem provides thus that
$$
[x^n]M(x)
\equivalent{n\to +\infty}\frac {2(e-1)}{1+4\pi ^2}\left( \frac 1{2\pi }\right) ^n \epsilon _n
$$
where
$$
\epsilon _n=
\left\{
\begin{array}{l}
1{\rm ~~if~}n{\rm ~is~even}\\
2\pi{\rm ~~if~}n{\rm ~is~odd}.
\end{array}
\right.
$$

\item Second case of mixing: $A=A'10^a$ and $B=0^b1B'$.

Because of Proposition~\ref{melange}, one has to compute $[x^n]M^{a,b}(x)$
with
$$
M^{a,b}(x):= S(1)\frac{c_{a+b}}{c_ac_b}P_{a+b}(x) 
+\frac 1{S(x)}\cdot\frac 1{1-x}\Big[ S(1)P_a(x)P_b(x)-S(x)\Big] ,
$$
where $ P_{a+b}$ is an entire function.
In this last formula, the brackets contain an entire function that vanishes at $1$ so that the dominant singularities
of $M^{a,b}$ are again those of $S^{-1}$, namely $\pm2i\pi$.
The expansion of $M^{a,b}(x)$ at $2i\pi$ writes thus
$$
M^{a,b}(x)\ \equivalent{2i\pi}\frac{-S(1)P_a(2i\pi )P_b(2i\pi )}{1-2i\pi}\cdot\frac{1}{1-\frac{x}{2i\pi}}
$$
which implies, by singularity analysis, that
$$
[x^n]M^{a,b}(x)\equivalent{n\to+\infty}2\Re\left(\frac{1-e}{1-2i\pi}\cdot\frac{P_a(2i\pi )P_b(2i\pi )}{(2i\pi )^n}\right) .
$$
Besides, the remainder of the exponential series satisfies
\begin{equation}
\label{resteSerieExp}
\sum _{n\geqslant a}\frac{x^n}{n!}=\frac{x^a}{a!}\left( 1+\frac xa+\mathcal{O}(\frac 1a)\right)
\end{equation}
when $a$ tends to infinity.
Consequently, by Formula~(\ref{defPa}), $P_a(2i\pi )$ tends to $2i\pi$ as $a$ tends to infinity
so that one gets a positive constant $C_1$ that does not depend on $a$ and $b$ such that for any
$n\geqslant 1$,
$$
\left| \psi(n, A,B)\right|\leqslant\frac{C_1}{(2\pi )^n}.
$$

\item Third case of mixing: $A=0^a$ and $B=0^b$.

This time, one has to compute $[x^n]M^{a,b}(x)$
with
$$
M^{a,b}(x) := S(1)\displaystyle\frac{1}{r_ar_b}\sum _{n\geqslant 1}r_{a+b+n-1}x^n + \displaystyle U(x)
\left[ \frac{S(1)R_a(x)R_b(x)}{r_ar_bx^{a+b-2}}-S(x)\right]
$$
the first term being an entire function.
Here again, the dominant singularities of $M^{a,b}$ are located at $\pm 2i\pi$ and
$$
M^{a,b}(x)\ \equivalent{2i\pi}\frac{-S(1)R_a(2i\pi )R_b(2i\pi )}{(1-2i\pi )r_ar_b(2i\pi )^{a+b-2}}\cdot\frac{1}{1-\frac{x}{2i\pi}}
$$
which implies, by singularity analysis, that
$$
\psi(n-1,A,B) \equivalent{n\to+\infty}
2\Re\left(\frac{1-e}{1-2i\pi}\cdot\frac{R_a(2i\pi )R_b(2i\pi )}{r_ar_b(2i\pi )^{a+b-2}}\frac{1}{(2i\pi )^n}\right) .
$$
Once more, because of~(\ref{resteSerieExp}), this implies that there is a positive constant $C_2$ independent of $a$ and $b$ and such that for any $n\geqslant 1$,
$$
\left| \psi(n,A,B)\right|\leqslant\frac{C_2}{(2\pi )^n}.
$$

\item and \itshape{\textbf{v)}}: both remaining cases of mixing that respectively correspond to words of the form $A=A'10^a$, $B=0^b$
and $A=0^a$, $B=0^b1B'$ are of the same vein and lead to similar results.\QED
\end{enumerate}

\section{Height and saturation level of suffix tries}
\label{sec:suffixtrie}

In this section, we consider a suffix trie process $(\rond{T}_n)_n$ associated with an infinite random word generated by an infinite comb. A precise definition of tries and suffix tries is given in section \ref{ssec:height}. We are interested in the height and the saturation level of such a suffix trie. 

Our method to study these two parameters uses a duality property \`a la Pittel developed in Section \ref{ssec:duality}, together with a careful and explicit calculation of the generating function of the second occurrence of a word (in Section \ref{ssec:generating}) which can be achieved for any infinite comb. These calculations are not so intricate because they are strongly related to the mixing coefficient and the mixing properties detailed in Section~\ref{sec:mixing}.

More specifically, we look at our two favourite examples, the logarithmic comb and the factorial comb.
We prove in Section \ref{ssec:heightLog} that the height of the first one is not logarithmic but polynomial  and in Section \ref{ssec:saturation} that the saturation level of the second one is not logarithmic either but negligibly smaller. Remark that despite the very particular form of the comb in the wide family of variable length Markov models, the comb sources provide a spectrum of asymptotic behaviours for the suffix tries.

\subsection{Suffix tries}
\label{ssec:height}
\begin{figure}[h]
 \begin{center}
 \includegraphics[height=5.5cm]{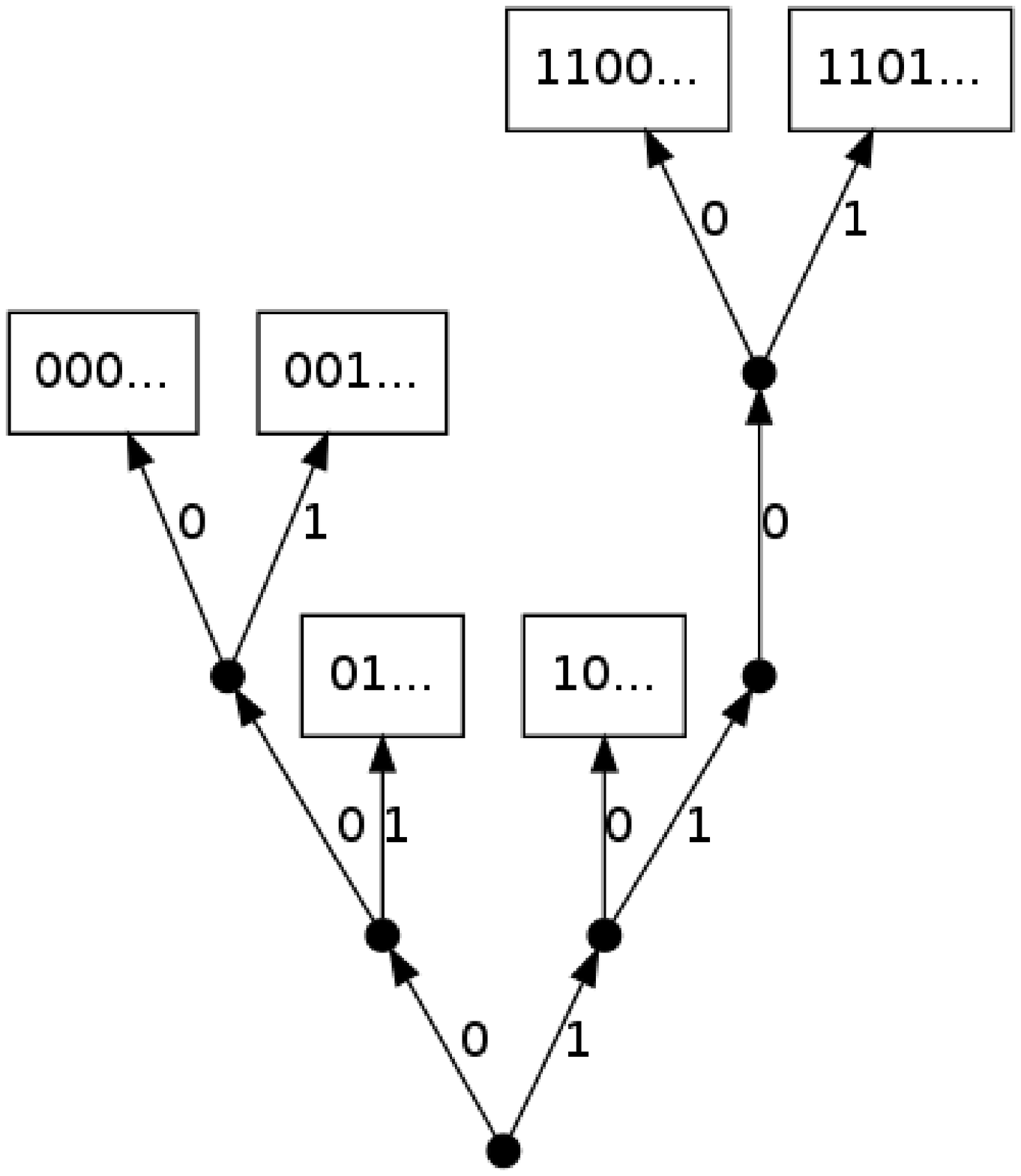} $\longrightarrow$ \includegraphics[height=7cm]{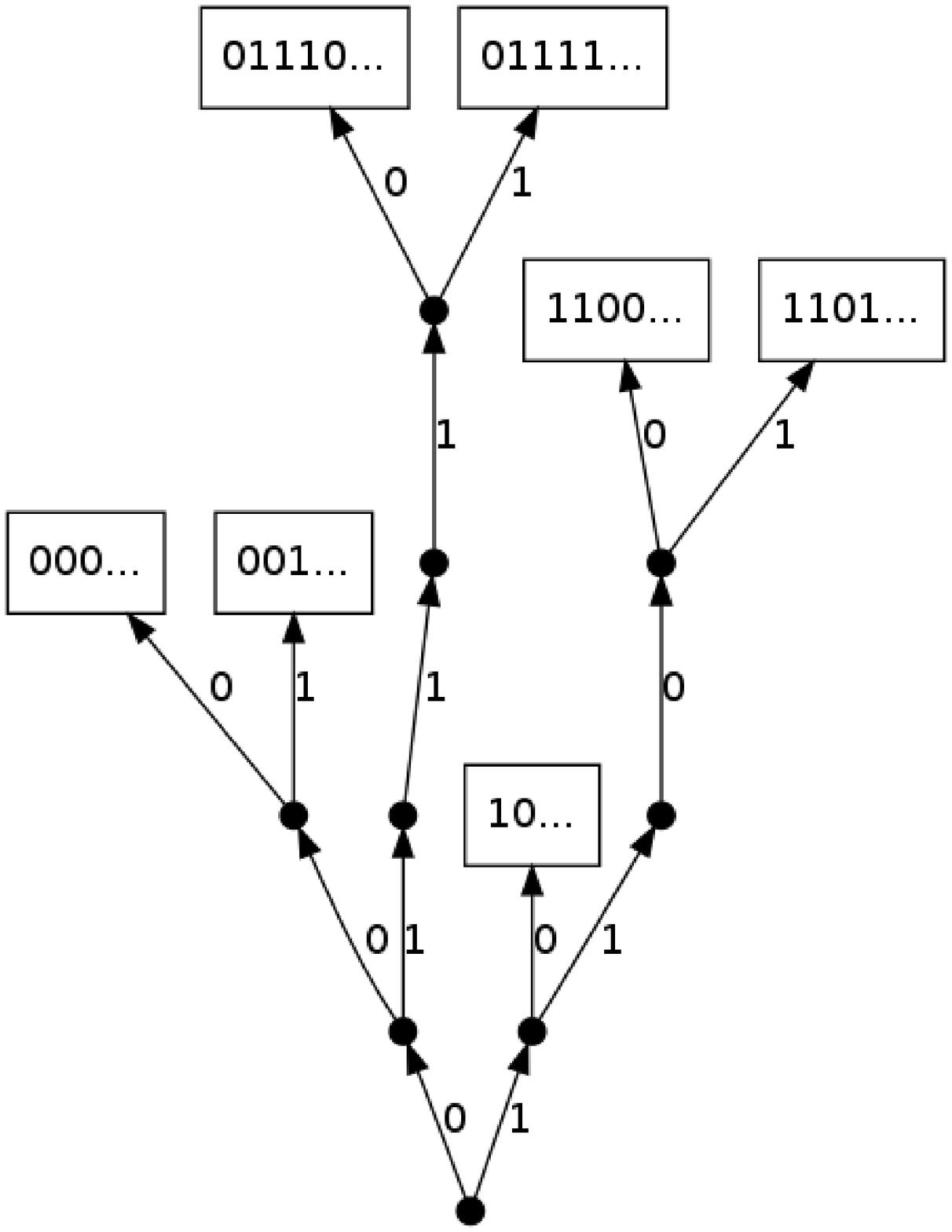} 
 \end{center}
\caption{\label{1figtrieconstruction}Last steps of the construction of a trie built from the set $(000\ldots,10\ldots, 1101\ldots, 001,\ldots, 01110\ldots , 1100\ldots, 01111\ldots)$.} 
\end{figure}

Let $(\rond{Y}_n)_{n\geq 1}$ be an increasing sequence of sets. Each set $\rond{Y}_n$ contains exactly $n$ infinite words. A \textbf{trie process} $(\rond{T}_n)_{n\geq 1}$ is a planar tree increasing process associated with $(\rond{Y}_n)_{n\geqslant 1}$. The trie $\rond{T}_n$ contains the words of $\rond{Y}_n$ in its leaves.
It is obtained by a sequential construction, inserting the words of $\rond{Y}_n$ successively. 
At the beginning, $\rond T_1$ is the tree containing the root and the leaf $0\dots$ (resp. the leaf $1\dots$) if the word in $\rond Y_1$ begins with $0$ (resp. with $1$). For $n\geq 2$, knowing the tree $\rond{T}_{n-1}$, the $n$-th word $m$ is inserted as follows. We go through the tree along the branch whose nodes are encoded by the successive prefixes of $m$; when the branch ends, if an internal node is reached, then the word is inserted at the free leaf, else we make the branch grow   comparing the next letters of both words until they can be inserted in two different leaves. As one can clearly see on Figure~\ref{1figtrieconstruction} a trie is not a complete tree and the insertion of a word can make a branch grow by more than one level. Notice that an internal node exists within the trie if there are at least two words in the set starting by the prefix associated to this node. This indicates why the \emph{second} occurrence of a word is prominent.

\medskip
 \begin{figure}[h]
 \begin{center}
 \psfrag{a}{\tiny{$10010\ldots$}}\psfrag{b}{\tiny{$00\ldots$}}\psfrag{c}{\tiny{$010\ldots$}}\psfrag{d}{\tiny{$101\ldots$}}\psfrag{e}{\tiny{$011\ldots$}}\psfrag{f}{\tiny{$11\ldots$}}\psfrag{g}{\tiny{$10011\ldots$}}\psfrag{00}{\tiny{$0$}}\psfrag{11}{\tiny{$1$}}
 \includegraphics[height=8.5cm]{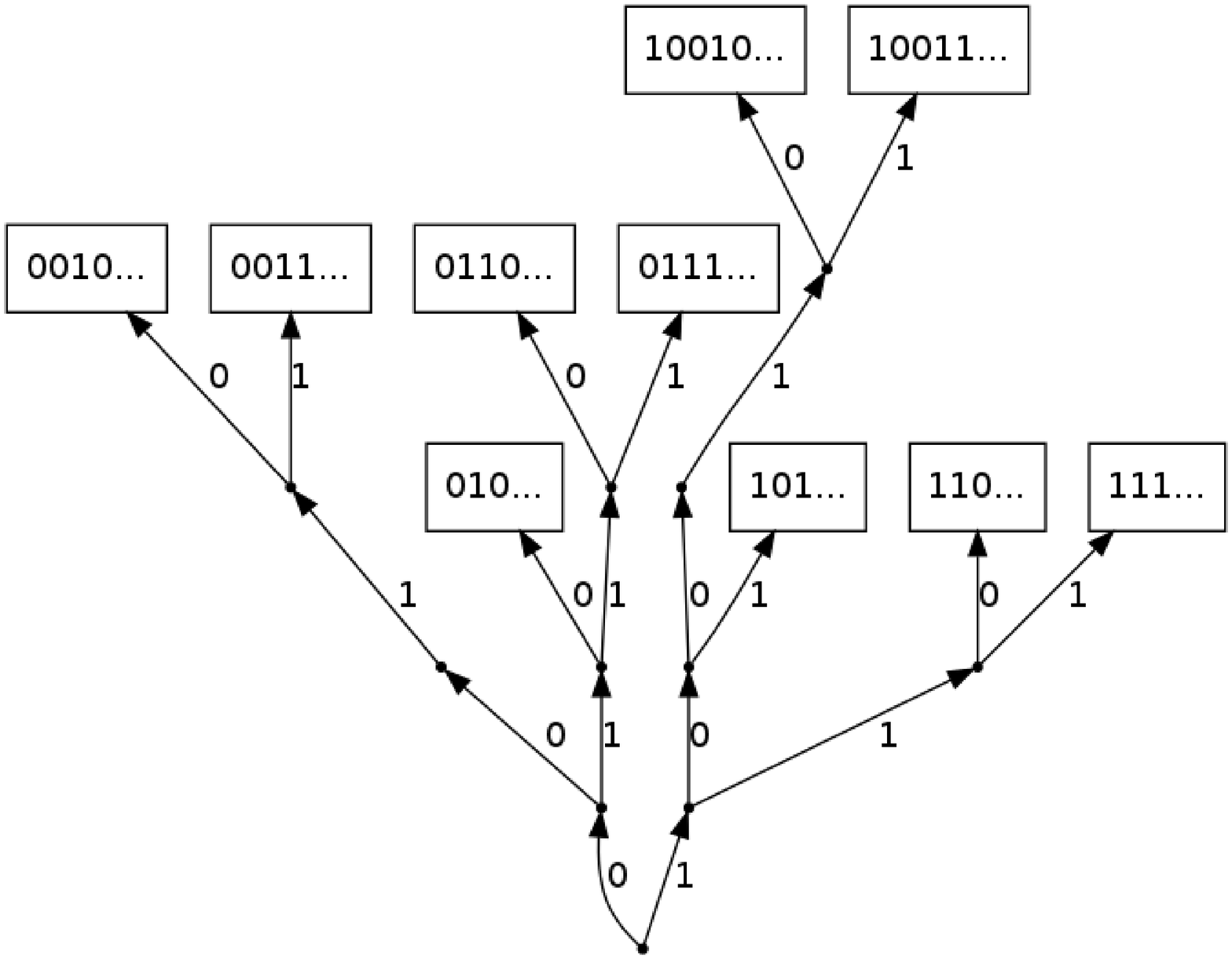}
 \end{center}
\caption{\label{1figtriesuffconstruction}Suffix trie $\rond T_{10}$ associated with the word $1001011001110\ldots$. Here, $H_{10}=4$ and $\ell_{10}=2$.}
 \end{figure}
\medskip

Let $m:=a_1a_2a_3 \ldots$ be an infinite word on $\rond{A}=\{0,1\}$. The \textbf{suffix trie} $\rond{T}_n$ (with $n$ leaves) associated with $m$, is the trie built from the set of the $n$-{th} first suffixes of $m$, that is 
\[
\rond{Y}_n = \{m,a_2a_3\ldots, a_3a_4\ldots, \ldots,a_na_{n+1}\ldots\}.\]
For a given trie  $\rond{T}_n$, we are mainly interested in  the \emph{height} $H_n$ which is the maximal depth of an internal node of $\rond{T}_n$ and the \emph{saturation level} $\ell_n$ which is the maximal depth up to which all the internal nodes are present in $\rond{T}_n$. Formally, if $\partial \rond T_n$ denotes the set of leaves of $\rond T_n$,
\begin{eqnarray*}
H_n &=&\max_{u \in \rond{T}_n\setminus \partial \rond{T}_n} \big\{ |u| \big\} \\
\ell_n&=&\max \big\{j \in \N | \ \# \{u \in  \rond{T}_n\setminus \partial \rond{T}_n, |u|=j\}=2^j\big\}.
\end{eqnarray*}
See Figure \ref{1figtriesuffconstruction} for an example.
\subsection{Duality}
\label{ssec:duality}

Let $(U_n)_{n\geqslant 1}$ be an infinite random word generated by some infinite comb and $(\rond{T}_n)_{n\geqslant 1}$ be the associated suffix trie process. We denote by $\rond R$ the set of right-infinite words.  Besides, we define hereunder two random variables having a key role in the proof of Theorem~\ref{hauteur} and Theorem~\ref{saturation}. This method goes back to Pittel \cite{Pittel}.

\begin{figure}[h]
\begin{center}
\psfrag{0}{\tiny{$0$}}\psfrag{1}{\tiny{$1$}}\psfrag{s}{\small{$s$}}\psfrag{t}{\small{$\hat{s}$}}
\includegraphics[height=9cm]{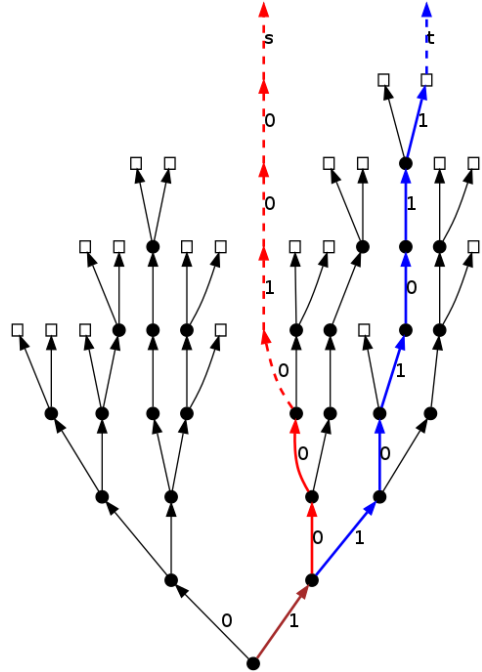}
\end{center}
\caption{\label{4figXnetTk}Example of suffix trie with $n=20$ words. The saturation level is reached for any sequence having $1000$ as prefix (in red); $\ell_{20}=X_{20}(s)=3$ and thus $T_3(s) \leqslant 20$. The height (related to the maximum of $X_{20}$) is realized for any sequence of the form $110101\dots$ (in blue) and $H_{20}=6$. Remark that the shortest branch has length $4$ whereas the saturation level $\ell_n$ is equal to $3$.}
\end{figure}

Let $s\in\rond R$ be a deterministic infinite sequence and $s ^{(k)}$ its prefix of length $k$. For $n\geq 1$,
\begin{eqnarray*}
\label{defXn}
X_{n}(s)&:=& \left \{
\begin{array}{l}
0 \ \mbox{if}\ s^{(1)}\ \mbox{is not in}\ \rond{T}_{n}\\
\max \{k\geqslant 1 \ |\ \mbox{the word}\ s ^{(k)} \mbox{ is already in }
\rond{T}_{n}\setminus \partial \rond{T}_n\}, \\
\end{array}
\right.\\
T_{k}(s)&:=& \min\{n\geqslant 1\ |\ X_{n}(s)=k\},
\end{eqnarray*}
where ``$s ^{(k)}$ is in $\rond{T}_{n}\setminus \partial \rond{T}_n$'' stands for: there exists an internal node $v$ in $\rond{T}_{n}$ such that $s ^{(k)}$ encodes $v$. For any $k\geqslant 1$, $T_{k}(s)$ denotes the number of leaves of the first tree ``containing'' $s^{(k)}$. See Figure \ref{4figXnetTk} for an example.
Thus, the saturation level $\ell_n$ and the height $H_n$ can be described using $X_n(s)$:
\begin{equation}
\label{defi-var}\ell_n=\min_{s\in \rond R} X_n(s) \quad \mbox{and}\quad H_n=\max_{s\in \rond R} X_n(s).\end{equation}
Moreover, $X_n(s)$ and $T_k(s)$ are in duality in the following sense: for all positive integers $k$ and $n$, one has the equality of the events
\begin{equation}\label{duality}
\{ X_{n}(s) \geqslant k \} = \{ T_{k}(s) \leqslant n \}.
\end{equation}
The random variable $T_{k}(s)$ (if $k\geqslant 2$) also represents the waiting time of the second occurrence of the deterministic word $s^{(k)}$ in the random sequence $(U_n)_{n\geq 1}$, \emph{i.e.}\@ one has to wait $T_{k}(s)$ for the source to create a prefix containing exactly two occurrences of $s^{(k)}$. More precisely, for $k \geqslant 2$, $T_{k}(s)$ can be rewritten as
\begin{align*}
T_k(s)=\min\Bigl\{n\geqslant 1\ \big| U_{n}U_{n+1}\ldots U_{n+k-1}&=s^{(k)} \mbox{  and } \exists ! j<n \mbox{ such that } \\
U_{j}U_{j+1}\ldots U_{j+k-1}&=s^{(k)}\Bigr\}.
\end{align*}
Notice that $T_k(s)$ denotes the \emph{beginning} of the second occurrence of $s^{(k)}$ whereas in \cite{ccpp}, $\tau^{(2)}\left(s^{(k)}\right)$ denotes the \emph{end} of the second occurrence of $s^{(k)}$, so that 
\begin{equation}
\label{lien-tau-T}
\tau^{(2)}\left(s^{(k)}\right)=T_k(s)+k.
\end{equation}
More generally, in \cite{ccpp} , for any $r\geqslant 1$,  the random return times $\tau^{(r)}(w)$ is defined as the end of the $r$-{th} occurrence of $w$ in the sequence $(U_n)_{n\geqslant 1}$ and the generating function of the $\tau^{(r)}$ is calculated. We go over these calculations in the sequel.
\subsection{Return time generating functions}
\label{ssec:generating}
\begin{Prop}
\label{Phi2}
Let $k\geqslant 1$. Let also $w=10^{k-1}$ and $\tau^{(2)}(w)$ be the end of the second occurrence of $w$ in a sequence generated by a comb defined by $(c_n)_{n\geqslant 0}$. Let $S$ and $U$ be the ordinary generating functions defined in Section \ref{ssec:mixingcomb}. The probability generating function of $\tau^{(2)}(w)$ is 
\[
 \Phi^{(2)}_w(x) = \frac{c_{k-1}^2x^{2k-1}\big(U(x)-1\big)}{S(1)(1-x)\big[ 1+c_{k-1}x^{k-1}(U(x)-1)\big] ^2}.
\]
Furthermore, as soon as $\sum_{n \geqslant 1}n^2c_n <\infty$, the random variable $\tau^{(2)}(w)$
is square-integrable and
\begin{equation}
\label{esp-var}
\g E(\tau^{(2)}(w))=\frac{2S(1)}{c_{k-1}}+o\left(\frac{1}{c_{k-1}}\right), \quad \Var(\tau^{(2)}(w))=\frac{2S(1)^2}{c_{k-1}^2}+o\left(\frac{1}{c_{k-1}^2}\right).
\end{equation}
\end{Prop}

\pff
For any $r\geqslant 1
$, let $\tau^{(r)}(w)$ denote the end of the $r$-th occurrence of $w$
in a random sequence generated by a comb and $\Phi_w^{(r)}$ its probability generating function. The reversed word of $c=\alpha _1\dots\alpha _N$ will be denoted by the overline $\overline c:=\alpha _N\dots\alpha _1$

We use a result of~\cite{ccpp} that computes these generating functions in terms of stationary
probabilities $q_c^{(n)}$. These probabilities measure the occurrence of a finite word after $n$ steps, conditioned
to start from the word $\overline{c}$.
More precisely,  for any finite words $u$ and $\overline{c}$ and for any
$n\geqslant 0$, let
\[q_c^{(n)}(u):=\pi\left(U_{n-|u|+|c|+1}\ldots U_{n+|c|}=u\big| U_{1}\ldots U_{|c|}=\overline c\right).\]
It is shown in \cite{ccpp} that, for $|x|<1$,
$$
\Phi_w^{(1)}(x)=\frac{x^k\pi(w)}{(1-x)S_w(x)}
$$
and for $r\geq 1$,
\[\Phi_w^{(r)}(x)=\Phi_w^{(1)}(x)\left(1-\frac{1}{S_w(x)}\right)^{r-1}
\]
where
\begin{eqnarray*}
S_w(x)&:=& C_w(x) + \sum_{n=k}^{\infty}q_{\petitlpref(w)}^{(n)}(w)x^n,\\
C_w(x)&:=& 1+\sum_{j=1}^{k-1}\ind{w_{j+1}\ldots w_{k}=w_{1}\ldots w_{k-j}}q_{\petitlpref(w)}^{(j)}\left(w_{k-j+1} \ldots w_k\right)x^j.
\end{eqnarray*}

In the particular case when $w=10^{k-1}$, then $\lpref(w)=\overline{w}=0^{k-1}1$ and $\pi(w)=\frac{c_{k-1}}{S(1)}$.
Moreover, Definition (\ref{mixingCoeff}) of the mixing coefficient and Proposition~\ref{melange}~\emph{\textbf{i)}} imply successively that
\begin{align*}
  q_{\petitlpref(w)}^{(n)}(w) & = \pi\Big(U_{n-k-|w|+1}\ldots U_{n+k}=w\Big| U_{1}\ldots U_k=\lpref(w)\Big) \\
   & = \pi(w) \Big( \psi\big( n-k,\lpref (w),w\big) +1\Big) \\
   & = \pi(w) S(1)u_{n-k+1}\\
   &= c_{k-1}u_{n-k+1},
   \end{align*}
This relation makes more explicit the link between return times and mixing.   
This leads to
 \[\sum_{n\geqslant k}q_{\petitlpref(w)}^{(n)}(w)x^n=c_{k-1}x^{k-1}\sum_{n\geqslant 1}u_nx^n=c_{k-1}x^{k-1}\big( U(x)-1\big) .\]
Furthermore, there is no auto-correlation structure inside $w$ so that $C_w(x)=1$ and
 \[S_w(x)=1+c_{k-1}x^{k-1}\big( U(x)-1\big) .\]
This entails
 \[\Phi^{(1)}_w(x)=\frac{c_{k-1}x^k}{S(1)(1-x)\left[ 1+c_{k-1}x^{k-1}\big( U(x)-1\big)\right]}\]
and
\begin{eqnarray*}
 \Phi^{(2)}_w(x) & = & \Phi^{(1)}_w(x)\left(1-\frac{1}{S_w(x)}\right) \\
  & = & \frac{c_{k-1}^2x^{2k-1}\big( U(x)-1\big)}{S(1)(1-x)\left[ 1+c_{k-1}x^{k-1}\big( U(x)-1\big)\right] ^2}
\end{eqnarray*}
which is the announced result.
The assumption
\[\sum_{n \geqslant 1}n^2c_n <\infty\]
makes $U$ twice differentiable and elementary calculations lead to 
\begin{eqnarray*}
&& (\Phi_w^{(1)})'(1)=\frac{S(1)}{c_{k-1}}-S(1)+1+\frac{S'(1)}{S(1)}, \qquad
(\Phi_w^{(2)})'(1)=(\Phi_w^{(1)})'(1)+\frac{S(1)}{c_{k-1}},\\
&& (\Phi_w^{(1)})''(1)=\frac{2S(1)^2}{c_{k-1}^2}+o\left(\frac{1}{c_{k-1}^2}\right)\quad \mbox{and} \quad (\Phi_w^{(2)})''(1)=\frac{6S(1)^2}{c_{k-1}^2}+o\left(\frac{1}{c_{k-1}^2}\right),\\
\end{eqnarray*}
and finally to (\ref{esp-var}).
\QED
\subsection{Logarithmic comb and factorial comb}
\label{ssec:results}
Let $h_+$ and $h_-$ be the constants in $[0,+\infty ]$ defined by
\begin{equation}
\label{4defh+} h_+ := \lim_{n\to +\infty}\frac{1}{n} \max\Bigl\{\ln\Bigl(\frac{1}{\pi\left(w\right)}\Bigr)\ \Bigr\}
\mbox{~and~~}h_- := \lim_{n\to +\infty}\frac{1}{n} \min\Bigl\{\ln\Bigl(\frac{1}{\pi\left(w\right)}\Bigr)\ \Bigr\},
\end{equation}
where the maximum and the minimum range over the words $w$ of length $n$ with $\pi\left({w}\right)>0$. In their papers, Pittel \cite{Pittel} and Szpankowski \cite{Szpan} only deal with the cases $h_+<+\infty$ and $h_->0$, which amounts to saying that the probability of any word is exponentially decreasing with its length. Here, we focus on our two examples for which these assumptions are not fulfilled. More precisely, for the logarithmic infinite comb, (\ref{probas}) implies that $\pi(10^n)$ is of order $n^{-4}$, so that
$$
h_- \leqslant \lim_{n\to +\infty}\frac{1}{n} \ln\Bigl(\frac{1}{\pi\left(10^{n-1}\right)}\Bigr) = 4\lim_{n\to +\infty}\frac{\ln n}{n} =0.
$$
Besides, for the factorial infinite comb, $\pi(10^n)$ is of order $\frac 1{(n+1)!}$ so that

$$
h_+ \geqslant \lim_{n\to +\infty}\frac{1}{n} \ln\Bigl(\frac{1}{\pi\left(10^{n-1}\right)}\Bigr) =  \lim_{n\to +\infty}\frac{n!}{n} = +\infty .
$$
For these two models, the asymptotic behaviour of the lengths of the branches is not always logarithmic, as can be seen in the two following theorems, shown in Sections~\ref{ssec:heightLog} and~\ref{ssec:saturation}.
\begin{Th}[Height of the logarithmic infinite comb] 
\label{hauteur} 
Let $\rond{T}_n$ be the suffix trie built from the $n$ first suffixes of a sequence generated by a logarithmic infinite comb. Then, the height $H_n$ of $\rond{T}_n$ satisfies
\[\forall \delta > 0, \hskip 1cm \frac{H_n}{n^{\frac 14 - \delta}}
~\limite{n\to \infty}{}~+\infty \quad \mbox{in probability.}
\]
\end{Th}


\begin{Th}[Saturation level 
of the factorial infinite comb] 
\label{saturation} 
Let $\rond{T}_n$ be the suffix trie built from the $n$ first suffixes of the sequence generated by a factorial infinite comb. Then, the saturation level $\ell_n$ of $\rond{T}_n$ satisfies:
for any $\delta >1$, almost surely, when $n$ tends to infinity,
$$
\ell _n\in o\left(\frac{\log n}{(\log\log n)^\delta}\right).
$$
\end{Th}
\begin{figure}[h]
 \begin{center}  
  \includegraphics[height=10cm]{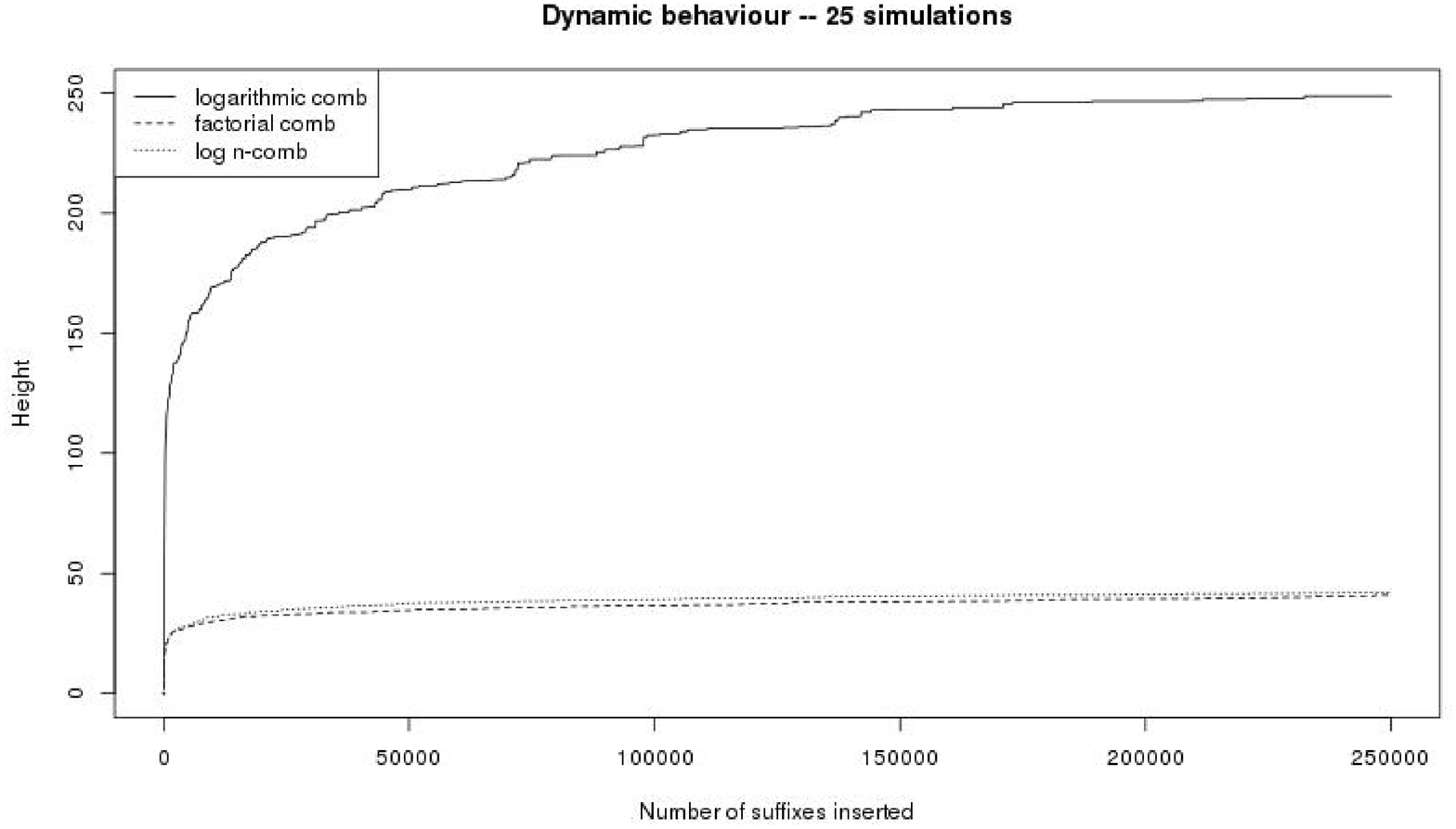}
  \includegraphics[height=10cm]{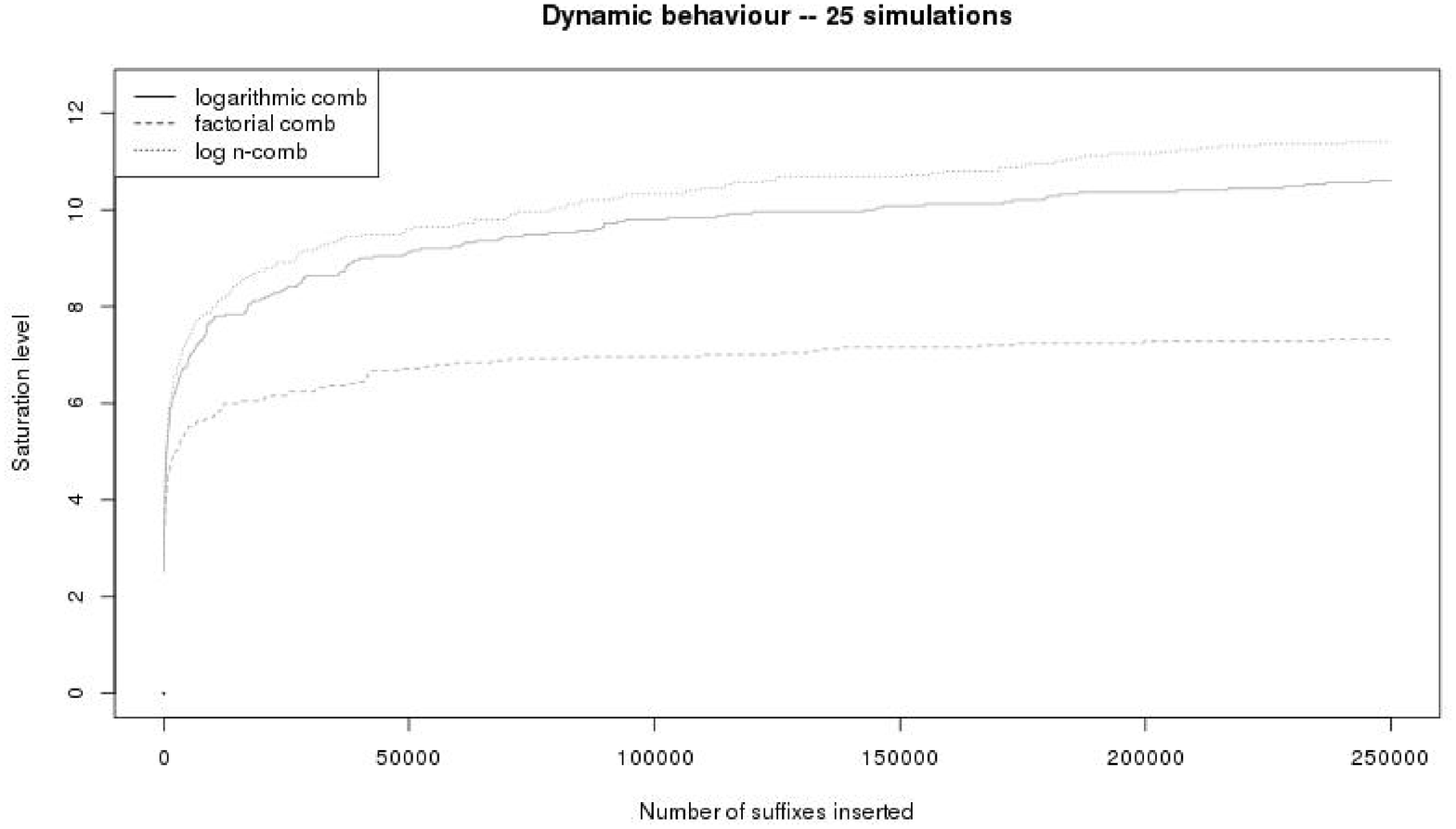}
 \end{center}
\caption{\label{dyn-peigne}
Respective heights and saturation levels for a logarithmic comb (plain lines), a factorial comb (long dashed lines) and a $\log n$-comb (short dashed lines).}
\end{figure}
The dynamic asymptotics of the height and of the saturation level can be visualized on Figure~\ref{dyn-peigne}.
The number $n$ of leaves of the suffix trie is put on the $x$-axis while heights or saturation levels of tries are put on the $y$-axis.
Plain lines represent a logarithmic comb while long dashed lines are those of a factorial comb (mean values of $25$ simulations).

Short dashed lines represent a third infinite comb defined by the data $c_n=\frac13\prod_{k=1}^{n-1}\left(\frac{1}{3}+\frac{1}{(1+k)^2}\right)$ for $n\geqslant 1$.
Such a process has a uniform exponential mixing, a finite $h_+$ and a positive $h_-$ as can be elementarily checked. As a matter of consequence, it satisfies all assumptions of Pittel \cite{Pittel} and Szpankowski \cite{Szpan} implying that the height and the saturation level are both of order $\log n$.
Such assumptions will always be fulfilled as soon as the data $(c_n)_n$ satisfy $\varlimsup _nc_n^{1/n}<1$; the proof of this result is left to the reader.

One can notice the height of the logarithmic comb that grows as a power of $n$.
The saturation level of the factorial comb, negligible with respect to $\log n$ is more difficult to highlight because of the very slow growth of logarithms.


\medskip
These asymptotic behaviours, all coming from the same model, the infinite comb, stress its surprising richness.


\subsection{Height for the logarithmic comb}
\label{ssec:heightLog}
In this subsection, we prove Theorem \ref{hauteur}. 

Consider the right-infinite sequence $s=10^{\infty}$. Then, $T_k(s)$ is the second occurrence time of $w=10^{k-1}$. It is a nondecreasing (random) function of $k$. Moreover, $X_n(s)$ is the maximum of all $k$ such that $s^{(k)}\in\rond T_n$. It is nondecreasing in~$n$. So, by definition of $X_n(s)$ and $T_k(s)$, the duality can be written
\begin{equation}\label{dualitybis}
\forall n, \forall \omega, \exists k_n, \quad k_n\leqslant X_n(s)<k_n+1 \mbox{~~and~~}T_{k_n}(s)\leqslant n < T_{k_n+1}(s). 
\end{equation}

{\bf 
Claim}: 
\begin{equation}\label{xninfini}
\lim_{n\rightarrow +\infty} X_n(s) = + \infty \quad \mbox{a.s.}
\end{equation}
Indeed, if $X_n(s)$ were bounded above, by $K$ say, then take $w = 10^{K}$ and consider $T_{K+1}(s)$ which is the time of the second occurrence of $10^{K}$.
The choice of the $c_n$ in the definition of the logarithmic comb implies the convergence of the series $\sum _nn^2c_n$. Thus (\ref{esp-var}) holds and $\g E[T_{K+1}(s)]<\infty$ so that $T_{K+1}(s)$ is almost surely finite. This means that for $n>T_{K+1}(s)$, the word $10^{K}$ has been seen twice, leading to $X_n(s)\geqslant K+1$ which is a contradiction.

\medskip
We make use of the following lemma that is proven hereunder.
\begin{Lem}\label{Tk} For $s=10^{\infty}$,
\[\forall \eta > 0, \quad \frac{T_k(s)}{k^{4+\eta}}~\limite{k \to \infty}{}~0{\sl ~~in~probability}, 
\]
and
\begin{equation}\label{Tksurk}
\forall \eta > \frac 12, \quad \frac{T_k(s)}{k^{4+\eta}} \ \ \smash{\mathop{\longrightarrow}\limits _{k\to \infty}} \ \ 0 \quad \mbox{a.s.}
\end{equation}
\end{Lem}

With notations~(\ref{dualitybis}), because of~(\ref{xninfini}), the sequence $(k_n)$ tends to infinity, so that $\left(T_{k_n}(s)\right)$ is a subsequence of $\left(T_k(s)\right)$. Thus, (\ref{Tksurk}) implies that
\[
\forall \eta > \frac 12, \quad \frac{T_{k_n}}{k_n^{4+\eta}} \ \ \smash{\mathop{\longrightarrow}\limits _{n\to \infty}} \ \ 0 \quad \mbox{a.s.} \quad \mbox{and} \quad  \forall \eta > 0, \quad \frac{T_{k_n}}{k_n^{4+\eta}}~\limite{k\to \infty}{P}~0.
\]
Using duality (\ref{dualitybis}) again leads to 
\[\forall \eta > 0, \quad \frac{X_n(s)}{n^{1/({4+\eta})}}~\limite{n \to \infty}{P}~+\infty.\]
In otherwords
$$
\forall \delta > 0, \quad \frac{X_n(s)}{n^{\frac 14 - \delta}}~\limite{n\to \infty}{P}~+\infty 
$$
so that, since the height of the suffix trie is larger than $X_n(s)$,
\[\forall \delta > 0, \quad \frac{H_n}{n^{\frac 14 - \delta}}~\limite{n\to \infty}{P}~+\infty.\]
This ends the proof of Theorem~\ref{hauteur}.
\QED

\bigskip
{\sc Proof of Lemma \ref{Tk}.}

Combining (\ref{lien-tau-T}) and (\ref{esp-var}) shows that
\begin{equation}\label{meanTk}
\g E(T_k(s)) = \g E(\tau^{(2)}(w))-k= \frac{19}{9} k^4 + o(k^4)
\end{equation}
and
\begin{equation}\label{varTk}
\Var(T_k(s)) = \Var(\tau^{(2)}(w))=\frac{361}{162} k^8 + o(k^8).
\end{equation}
For all $\eta >0$, write
$$
\frac{T_k(s)}{k^{4+\eta}} = \frac{T_k(s) -\g E(T_k(s))}{k^{4+\eta}} + \frac{\g E(T_k(s))}{k^{4+\eta}} .
$$
The deterministic part in the second-hand right term goes to $0$ with $k$ thanks to (\ref{meanTk}), so that we focus on the term $\displaystyle\frac{T_k(s) -\g E(T_k(s))}{k^{4+\eta}}$.
For any $\varepsilon >0$, because of Bienaym\'e-Tchebychev inequality,
$$
\g \PP\left(\frac{|T_k(s) -\g E[T_k(s)]|}{k^{4+\eta}} >\varepsilon\right) \leqslant \frac{Var(T_k(s))}{\varepsilon^2 k^{8+2\eta}}=\mathcal{O}\left(\frac 1{k^{2\eta}}\right).
$$
This shows the convergence in probability in Lemma~\ref{Tk}. Moreover, Borel-Cantelli Lemma ensures the almost sure convergence as soon as $\eta > \frac 12$. \QED

\begin{Rem}
Notice that our proof shows actually that the convergence to $+\infty$ in Theorem~\ref{hauteur} is valid a.s. (and not only in probability) as soon as $\delta > \frac 1{36}$. \end{Rem}

\subsection{Saturation level for the factorial comb}
\label{ssec:saturation}
In this subsection, we prove Theorem \ref{saturation}.

Consider the probabilized infinite \emph{factorial} comb defined in Section~\ref{sec:source}
by
$$
\forall n\in\g N,~c_n=\frac{1}{(n+1)!}.
$$
The proof hereunder shows actually that $\left( \frac{\ell _n\log\log n}{\log n}\right) _n$ is an almost surely bounded sequence, which implies the result.
Recall that $\rond R$ denotes the set of all right-infinite sequences.
By characterization of the saturation level as a function of $X_n$ (see (\ref{defi-var})),
$\g P\left(\ell _n\leqslant k\right) =\g P\left(\exists s\in\rond R,~X_n(s)\leqslant k\right)$ for all positive integers $n,k$.
Duality formula~(\ref{duality}) then provides
$$
\begin{array}{rl}
\g P\left(\ell _n\leqslant k\right)&=\g P\left(\exists s\in\rond R,~T_k(s)\geqslant n\right)\\ \\
&\geqslant\g P\left( T_k( \widetilde{s})\geqslant n\right)
\end{array}
$$
where $\widetilde{s}$ denotes any infinite word having $10^{k-1}$ as a prefix.
Markov inequality implies
\begin{equation}
\forall x\in ]0,1[,~\g P\left( \ell _n\geqslant k+1\right)\leqslant \g P\Big(\tau^{(2)}(10^{k-1})<n+k\Big)\leqslant\frac{\Phi ^{(2)}_{10^{k-1}}(x)}{x^{n+k}}
\end{equation}
where $\Phi ^{(2)}_{10^{k-1}}(x)$
 denotes as above the generating function of the rank of the final letter of the second occurrence of $10^{k-1}$ in the infinite random word $(U_n)_{n\geq 1}$. The simple form of the factorial comb leads to the explicit expression
$U(x)=\frac{x}{(1-x)(e^x-1)}$ and, after computation,
\begin{equation}
\label{phi2Factoriel}
\Phi ^{(2)}_{10^{k-1}}(x)=
\frac{e^x-1}{e-1}\cdot
\frac{x^{2k-1}\big( 1-e^x(1-x)\big)}
{\Big[ k!\left( e^x-1\right)\left( 1-x\right)+x^{k-1}\big( 1-e^x(1-x)\big)\Big] ^2}.
\end{equation}
In particular, applying Formula~(\ref{phi2Factoriel}) with $n=(k-1)!$ and $x=1-\frac 1{(k-1)!}$ implies that
for any $k\geqslant 1$,
$$
\g P(\ell _{(k-1)!}\geqslant k+1)\leqslant
\frac{\left( 1-\frac{1}{(k-1)!}\right) ^{2k-1}}{\Big[ k!(e-1)\frac{1}{(k-1)!}\Big] ^2}
\cdot\frac {1}{\left( 1-\frac{1}{(k-1)!}\right) ^{(k-1)!+k}}.
$$
Consequently, $\g P\left(\ell _{(k-1)!}\geqslant k+1\right)=\mathcal{O}(k^{-2})$ is the general term of a convergent series.
Thanks to Borel-Cantelli Lemma, one gets almost surely
\[
\varlimsup _{n\to+\infty}\frac{\ell _{n!}}{n}\leqslant 1.
\]
Let $\Gamma ^{-1}$ denote the inverse of Euler's Gamma function, defined and increasing on the real interval
$[2,+\infty [$.
If $n$ and $k$ are integers such that $(k+1)!\leqslant n\leqslant (k+2)!$, then
$$
\frac{\ell _n}{\Gamma ^{-1}(n)}\leqslant\frac{\ell _{(k+2)!}}{\Gamma ^{-1}((k+1)!)}=\frac{\ell _{(k+2)!}}{k+2},
$$
which implies that, almost surely,
$$
\varlimsup _{n\to\infty}\frac{\ell _n}{\Gamma ^{-1}(n)}\leqslant 1.
$$
Inverting Stirling Formula, namely
$$
\Gamma (x)=\sqrt{\frac{2\pi}{x}}e^{x\log x-x}\left( 1+\mathcal{O}\left( \frac 1x\right)\right)
$$
when $x$ goes to infinity, leads to the equivalent
$$
\Gamma ^{-1}(x)\equivalent {+\infty}\frac{\log x}{\log\log x},
$$
which implies the result.
\QED

\subsection*{Acknowledgements}
The authors are very grateful to Eng. Maxence Guesdon for providing simulations with great talent and an infinite patience.
They would like to thank also all people managing two very important tools for french mathematicians: first the Institut Henri Poincar\'e, where a large part of this work was done and second Mathrice which provides a large number of services.

\bibliographystyle{plain} 
\bibliography{vlmc}

\end{document}